# ОБОБЩЕННОЕ ТЕОРЕТИКО-ЧИСЛОВОЕ ПРЕОБРАЗОВАНИЕ

Мирослав Семотюк

Институт кибернетики имени В.М. Глушкова НАН Украины,
проспект Академика Глушкова, 40, Киев 03187, Украина
e-mail: semo@i.ua, факс: +38 044 5263348.

**Аннотация**

Обобщенное теоретико-числовое преобразование (ТЧП) сформулировано на основе теоремы о показательной функции, которая позволяет заменять операции по модулю выражения в целом операциями по модулю над показателем степени этой функции, что составляет фундаментальность этой теоремы для ТЧП, поскольку именно такая функция используется в ТЧП как весовая функция преобразования. На основании этой теоремы в работе сформулированы и доказаны все основные теоремы обобщенного ТЧП, их двойственность а также свойства весовых функций этого преобразования. Выбор основания этой функции, в качестве которого может быть выбрано любое число, в том числе и комплексное, определяет не только тот или иной вид преобразования, но и модуль самого преобразования. Это позволяет обобщить ряд известных ТЧП, таких как Мерсена, Гаусса и даже Фурье в виде единой теории дискретных преобразований.

**Ключевые слова:** обобщенное теоретико-числовое преобразование, изображение, оригинал, теорема, вычет, модуль, кольцо вычетов, свертка, весовая функция, корреляционная функция, ядро преобразования.

The generalized number-theoretic transformation (NPT) is formulated on the basis of the exponential function theorem, which allows us to replace operations modulo the expression as a whole by modulo operations on the exponent of this function, which makes this theorem fundamental for NPT, since it is such a function used in NPT as a weight conversion function. On the basis of this theorem, all the main theorems of the generalized NPT, their duality, as well as the properties of the weight functions of this transformation are formulated and proved. The choice of the basis of this function, as which any number can be chosen, including a complex one, determines not only one or another type of transformation, but also the module of the transformation itself. This allows us to generalize a number of well-known NPTs, such as Mersen, Gauss, and even Fourier, in the form of a unified theory of discrete transformations.

## 1. Общие сведения

Алгоритмы обработки сигналов в конечных математических структурах, аналогично спектральной обработке, составляют интересную область исследований. Их применение, вопреки сложившемуся мнению, несомненно, не ограничиваются лишь вычислением сверток. Они полезны при описании систем счисления, анализе кодовых последовательностей, операциях над полиномами, а также для решения ряда задач, где требуются точные значения результатов вычислений. Однако возможности теоретико-числовых преобразований (а именно так называют обработку в конечных математических структурах), первоначально казавшиеся весьма многообещающими непосредственно для обработки сигналов, оказались сильно ограниченными и интерес специалистов, занимающихся такой обработкой, постепенно угас.

Вместе с тем, несмотря на сравнительно большое количество публикаций по этому вопросу [1,2,3], достаточно полного описания теоретико-числовых преобразований с полными выкладками их теорем, как это сделано для Фурье - преобразований или родственных ему, в литературе нет. Многие публикации отображают лишь отдельные фрагменты преобразования, или предлагают читателю самому доказать теорему о свертке, или искать подходящие модули преобразования и т.п.

Таким образом, на основании выше изложенного возникает необходимость в создании некоторой обобщенной теории, специально построенной для теоретико-числовых преобразований, которая систематизировала бы ряд полученных за последнее время результатов.

В качестве основы этой теории, по мнению автора, может быть использована фундаментальная теорема теоретико-числового преобразования [4], доказательство которой приводится ниже. А в силу того, что дальнейшие результаты получены только на ее основе, то ссылки на литературные источники делаться не будут. Автор в данной работе также придерживается финитной точки зрения в математике, согласно которой все рассматриваемые объекты не выходят за пределы некоторого ограниченного множества, на котором суждения о них истинны.

## 2. Фундаментальная теорема теоретико-числового преобразования.

Сформулируем следующую теорему. Пусть алгебра вида $\mathcal{K}_m = <\mathbf{S}, +, \cdot, \mathbf{0}, \mathbf{1}>$ где ($\mathbf{0} \neq \mathbf{1}$) и $\mathbf{S} \in \mathbf{Z}$ - структура (решетка) имеющая $\sup \mathbf{S} = s^p - 1$, $\inf \mathbf{S} = 0$ представляет собой кольцо вычетов с единицей, в котором своими аргументами задана степенная зависимость $y = s^x$. Тогда для $\forall p \in \mathbf{N}, \forall x = \overline{0, N}$ и $\forall p \ll N$ существует число $M < \sup \mathbf{S}$, при котором справедливо следующее равенство в кольце вычетов $\mathcal{K}_m$.

$$s^{(x) \bmod p} \stackrel{\mathcal{K}_m}{=\!=} (s^x) \bmod M, \qquad (2.1)$$

где $\stackrel{\mathcal{K}_m}{=\!=}$ - обозначение равенства в кольце вычетов (имеется в виду финитная точка зрения в математике), которое не всегда совпадает с известным понятием «сравнение по модулю» в силу разных значений модуля в левой и правой частях выражения (2.1). Модуль $M$ при этом есть функция от переменной $p$ - $M = f(p)$. Другими словами, если в левой или правой части равенства есть операция по модулю, то независимо от нее результат еще раз ограничивается, как слева, так и справа равенства, модулем кольца $\mathcal{K}_m$.

Докажем это утверждение, полагая первоначально, что это число равно верхней грани структуры $\mathbf{S}$
$$M = s^p - 1. \qquad (2.2)$$
Тогда, с одной стороны,

$x = kp + (x) \bmod p$, где $k = \mathrm{int}(x / p)$ - целая часть от деления $x$ на $p$;

$(x) \bmod p$ - остаток от этого же деления.

Отсюда функция $y = s^x$ может быть представлена

$$s^x = s^{kp + (x) \bmod p}. \qquad (2.3)$$

С другой стороны, аналогично

$$s^x = c(s^p - 1) + (s^x) \bmod (s^p - 1). \qquad (2.4)$$

Приравнивая (2.3) и (2.4), имеем

$$s^{kp+(x) \bmod p} = c(s^p - 1) + (s^x) \bmod (s^p - 1)$$

или

$$s^{kp} s^{(x) \bmod p} = c(s^p - 1) + (s^x) \bmod (s^p - 1). \qquad (2.5)$$

Заметим, что равенство (2.1.1) с учетом (2.1.2) может быть достигнуто при следующих соотношениях

$$s^{kp} = 1, \; c(s^p - 1) = 0. \qquad (2.6)$$

Из формулы (2.1) следует, что функция $y = s^x$ в правой части этого выражения при изменении $x = \overline{0, N}$, ограничена сверху модулем $M$. Очевидно, чтобы сохранить корректность результатов, необходимо ограничить и левую часть выражения (2.1). Тогда (2.5) можно записать так

$$(s^{kp} s^{(x) \bmod p}) \bmod (s^p - 1) = c(s^p - 1) + (s^x) \bmod (s^p - 1).$$

На основании последнего, выражения (2.6) могут быть заданы сравнениями вида

$$s^{kp} = 1 \,|\, \bmod (s^p - 1), \qquad (2.7)$$

$$c(s^p - 1) = 0 \,|\, \bmod (s^p - 1), \qquad (2.8)$$

где вертикальная черта «|» - знак, разделяющий основное выражение и модуль сравнения.

Тождество (2.8) справедливо в силу линейности сравнения, а (2.7) можно переписать таким образом

$$s^{pk} - 1 = 0 \,|\, \bmod (s^p - 1). \qquad (2.9)$$

По теореме Безу левую часть последнего, как разность степеней, можно разложить на сомножители вида

$$s^{pk} - 1 = (s^p)^k - 1 = (s^p - 1) \cdot F(s^p)^{k-1}, \text{ где } F(s^p)^{k-1} \text{ - полином } (k-1)\text{-степени от } s^p.$$

Так как $s^x$ является целочисленной функцией, то $s^{pk} - 1$ делится на $s^p - 1$ без остатка и

$$s^p - 1 = 0 \,|\, \bmod (s^p - 1),$$

откуда следует справедливость тождества (2.7). Теперь (2.5) можно представить в виде:

$$1 \cdot s^{(x) \bmod p} = 0 + (s^x) \bmod (s^p - 1) \,|\, \bmod (s^p - 1)$$

или

$$s^{(x) \bmod p} = (s^x) \bmod (s^p - 1) \,|\, \bmod (s^p - 1).$$

Полагая, что модуль кольца $\mathcal{K}_m$ равен модулю сравнения, последнее можно записать в следующем виде

$$s^{(x) \bmod p} \stackrel{\mathcal{K}_m}{=\!=\!=} (s^x) \bmod (s^p - 1). \qquad (2.10)$$

Таким образом, найдено значение числа $M$, удовлетворяющее теореме (2.1), однако совпадающее с модулем кольца $\mathcal{K}_m$. Покажем, что это значение не является единственным. Действительно, умножив и разделив это число на $s-1$, будем иметь

$$M = \frac{s^p - 1}{s - 1}(s - 1).$$

$\frac{s^{p-1}}{s-1}$ - является известным выражением для суммы ряда геометрической прогрессии со знаменателем $s$ и свободным членом $a_0 = 1$. Тогда $M$ является составным числом.

$$M = s^p - 1 = (s - 1) \sum_{m=0}^{p-1} s^m. \qquad (2.11)$$

Рассмотрим теперь число $\sum_{m=0}^{p-1} s^m$. Из (2.11) ясно, что оно удовлетворяет следующему соотношению

$$s^p - 1 < \sum_{m=0}^{p-1} s^m < s^p.$$

Однако этому же неравенству удовлетворяет также число $M$

$$s^{p-1} < s^p - 1 < s^p.$$

Стало быть, эти числа имеют один и тот же порядок относительно $s^x$ и, следовательно, остатки от деления $s^x$ на эти числа совпадают. Действительно

$$(s^x) \bmod (s^p - 1) = s^x - c(s^p - 1), \text{ где } c = \operatorname{int} \frac{s^x}{s^p - 1},$$

или с учетом (2.11) $\quad (s^x) \bmod (s^p - 1) = s^x - c(s-1) \cdot \sum_{m=0}^{p-1} s^m$, где $c(s-1) = \operatorname{int} \dfrac{s^x}{\sum_{m=0}^{p-1} s^m}$.

Тогда $\quad (s^x) \bmod (s^p - 1) \stackrel{\mathcal{K}_m}{=\!=\!=} s^x \bmod (\sum_{m=0}^{p-1} s^m)$.

Таким образом, найдено еще одно значение числа $M$, которое не совпадает с модулем кольца $\mathcal{K}_m$

$$M = \sum_{m=0}^{p-1} s^m. \qquad (2.12)$$

Отсюда следует, что теорема (2.1) доказана. Отметим некоторые следствия, вытекающие из этой теоремы.

*Следствие* 1. Числа $s$ и $\sum_{m=0}^{p-1} s^m$ взаимно простые при $\forall s > 2$.

Действительно $\quad \sum_{m=0}^{p-1} s^m = s^0 + \sum_{m=1}^{p-1} s^m = s^0 + \sum_{m=0}^{p-2} s^m$.

Поскольку $s < \sum_{m=0}^{p-1} s^m$, то, учитывая последнее, будем иметь

$$\frac{\sum_{m=0}^{p-1} s^m}{s} = \frac{s^0}{s} + \sum_{m=0}^{p-2} s^m, \quad \text{здесь} \quad \sum_{m=0}^{p-2} s^m = \text{int} \frac{\sum_{m=0}^{p-1} s^m}{s}, \quad \text{а} \quad \frac{s^0}{s} = \frac{1}{s} < 1.$$

Тогда $\forall s > 2$ не является делителем нуля для выражения $\sum_{m=0}^{p-1} s^m$ и, следовательно, эти числа являются взаимно простыми.

*Следствие* 2. Две числовые последовательности вида $s^{(x) \bmod p}$ и $(s^x) \bmod (\sum_{m=0}^{p-1} s^m)$ полученные из функции $s^x$ путем изменения $x = \overline{0, N}$ в кольце вычетов $\mathcal{K}_m$ конгруэнтны вплоть до каждого члена при одном и том же значении $x$.

*Следствие* 3. Две числовые последовательности вида $s^{(x) \bmod p}$ и $(s^x) \bmod (\sum_{m=0}^{p-1} s^m)$ при $x = \overline{0, N}$ периодичны в кольце вычетов $\mathcal{K}_m$, имеют одинаковый период $p$ и одно и то же главное значение, находящееся в интервале $[0, p-1]$.

Отметим теперь еще одно важное свойство теоремы (2.1). Эта теорема позволяет заменять операции по модулю выражений в целом операциями по модулю над показателями степеней степенных зависимостей, входящих в эти выражения, что и составляет фундаментальность теоремы (2.1) для теоретико-числовых преобразований, а также является отправным моментом для доказательства многих основных положений таких преобразований.

### 3. Обобщенное теоретико-числовое преобразование (**S**-преобразование).

Полагая, что главное значение числовой степенной последовательности находится на закрытом интервале $[0, p-1] = [0, N-1]$, при этом $\sup \mathbf{S} = M$, где $M = \sum_{m=0}^{N-1} s^m$ - модуль кольца $\mathcal{K}_m$, определим формально следующее преобразование, заданное на структуре **S**

$$X(k) \stackrel{\mathcal{K}_m}{=\!=} \sum_{i=0}^{N-1} x(i) s^{-(ki) \bmod N}, \tag{3.1}$$

$$x(i) \stackrel{\mathcal{K}_m}{=\!=\!=} \frac{1}{N}\sum_{k=0}^{N-1} X(k) s^{(ki)\bmod N}, \qquad (3.2)$$

$$M = \sum_{m=0}^{N-1} s^m,$$

где $\mathcal{K}_m$ - кольцо вычетов по модулю $M$, $N$ - некоторое число из множества **N**,

$\stackrel{\mathcal{K}_m}{=\!=\!=}$ - означает равенство (сравнение) в кольце вычетов $\mathcal{K}_m$,

$x(i), X(k)$ - числовые последовательности, представляющие оригинал и изображение соответственно, $s$ - некоторое число, в общем случае комплексное, $i,k$ - номера (индексы) компонент последовательностей.

Выражение (3.1) представляет собой прямое преобразование, а (3.2) - обратное.

Покажем теперь, что эта пара выражений действительно представляет собой теоретико-числовое преобразование. Для чего, во избежание путаницы индексов в выражении (3.2), заменим переменную $i$ переменной $n$, а затем, подставив ее в (3.1), в результате будем иметь

$$x(n) \stackrel{\mathcal{K}_m}{=\!=\!=} \frac{1}{N}\sum_{k=0}^{N-1}\left(\sum_{i=0}^{N-1} x(i) s^{-(ki)\bmod N}\right) s^{(kn)\bmod N}.$$

Изменив порядок суммирования, получим

$$x(n) \stackrel{\mathcal{K}_m}{=\!=\!=} \frac{1}{N}\sum_{k=0}^{N-1} x(i)\left(\sum_{i=0}^{N-1} x(i) s^{(kn-ki)\bmod N}\right).$$

Далее $\qquad x(n) \stackrel{\mathcal{K}_m}{=\!=\!=} \frac{1}{N}\sum_{i=0}^{N-1} x(i)\left(\frac{1}{N}\sum_{k=0}^{N-1} s^{|(n-i)k|\bmod N}\right). \qquad (3.3)$

Рассмотрим теперь внутреннюю сумму

$$\sum_{k=0}^{N-1} s^{|(n-i)k|\bmod N}. \qquad (3.4)$$

Очевидно, что при $i=n$, значение этой суммы будет равно $N$. Для $i \neq n$ внутренняя сумма (3.4) не равна $N$ и не равна нулю, что необходимо для преобразования. Однако, может существовать сравнение вида

$$\sum_{k=0}^{N-1} s^{|(n-i)k|\bmod N} = 0 \,|\bmod M, \qquad (3.5)$$

которого достаточно для преобразования, поскольку вычисления будут выполняться в кольце вычетов $\mathcal{K}_m$. Положим, что $M = \sum_{m=0}^{N-1} s^m$. Тогда, на основании теоремы (3.1) имеем

$$\sum_{k=0}^{N-1} s^{|(n-i)k|\bmod N} \stackrel{\mathcal{K}_m}{=\!=\!=} \left(\sum_{k=0}^{N-1} s^{(n-i)k}\right)\bmod M,$$

а с учётом (3.5) будем иметь сравнение $\sum_{k=0}^{N-1} s^{|(n-i)k} = 0\,|\,\mathrm{mod}\,M$

или $$\sum_{k=0}^{N-1} (s^{n-i})^k = 0\,|\,\mathrm{mod}\,M.$$

На основании (3.1) можно записать $\dfrac{s^{(n-i)N}-1}{s^{n-i}-1} = 0\,|\,\mathrm{mod}\,M$

Или $$s^{(n-i)N} = 0\,|\,\mathrm{mod}\,M.$$

Тогда сравнение примет вид $$\dfrac{s^{(n-i)N}-1}{s^{n-i}-1} = 0\,|\,\mathrm{mod}\,M. \tag{3.6}$$

На основании теоремы (3.1) можно записать

$$(s^{(n-i)N})\,\mathrm{mod}\,M \xlongequal{\mathcal{X}_m} s^{[(n-i)N]\,\mathrm{mod}\,N} \xlongequal{\mathcal{X}_m} s^0, \qquad (s^{(n-i)N})\,\mathrm{mod}\,M \xlongequal{\mathcal{X}_m} 1.$$

Подставляя в (3.6), будем иметь $\quad 1 - 1 = 0\,|\,\mathrm{mod}\,M$.

Стало быть, сравнение (3.6) выполняется и, следовательно, выполняется сравнение (3.5).

Таким образом, установлено, что для всех $i = n$ внутренняя сумма (3.4) равна $N$, а при $i \neq n$ она равна нулю в кольце вычетов $\mathcal{X}_m$. Последнее, с учётом выше сказанного, а также порядком выполнения операций в кольце вычетов будет иметь вид следующей системы уравнений

$$\begin{cases} \sum_{m=0}^{N-1} s^{[(n-i)k]\,\mathrm{mod}\,N} \xlongequal{\mathcal{X}_m} N, \; i = n \\ \sum_{m=0}^{N-1} s^{[(n-i)k]\,\mathrm{mod}\,N} \xlongequal{\mathcal{X}_m} 0, \; i \neq n \end{cases}, \tag{3.7}$$

Отсюда следует, что и правая часть выражения (3.3) будет состоять из единственного, не равного нулю в кольце вычетов $\mathcal{X}_m$, члена $x(i)$ только в том случае, если $i = n$. Тогда в (3.3) последовательность $x(n)$ совпадает с последовательностью $x(i)$ при $i = n$, следовательно, выражения (3.1) и (3.2) представляют собой теоретико-числовое преобразование. Последнее выражение легко иллюстрируется в матричном виде.

Пусть размер преобразования $N = 3$, а модуль $M = \sum_{m=0}^{2} s^m$. В этом случае имеем две матрицы. Одну матрицу для прямого преобразования, другую для обратного преобразования.

$$\begin{bmatrix} 1 & 1 & 1 \\ 1 & s^1 & s^2 \\ 1 & s^2 & s^1 \end{bmatrix}, \begin{bmatrix} 1 & 1 & 1 \\ 1 & 1/s^1 & 1/s^2 \\ 1 & 1/s^2 & 1/s^1 \end{bmatrix}$$

Найдем обратные элементы для величин $1/s^1$ и $1/s^2$. Ими будут, соответственно, $s^2$ и $s^1$ в силу того, что $s^2 \cdot s^1 = 1 | \mod M$. Тогда матрицы примут вид

$$\begin{bmatrix} 1 & 1 & 1 \\ 1 & s^1 & s^2 \\ 1 & s^2 & s^1 \end{bmatrix}, \begin{bmatrix} 1 & 1 & 1 \\ 1 & s^2 & s^1 \\ 1 & s^1 & s^2 \end{bmatrix}$$

а их произведение в кольце вычетов $\mathcal{K}_m$ будет равно

$$\left[ \begin{bmatrix} 1 & 1 & 1 \\ 1 & s^1 & s^2 \\ 1 & s^2 & s^1 \end{bmatrix} \times \begin{bmatrix} 1 & 1 & 1 \\ 1 & s^2 & s^1 \\ 1 & s^1 & s^2 \end{bmatrix} \right] \mod (\sum_{m=0}^{2} s^m) = \begin{bmatrix} 3 & 0 & 0 \\ 0 & 3 & 0 \\ 0 & 0 & 3 \end{bmatrix}$$

в чем нетрудно убедиться, выполнив соответственно вычисления, принимая во внимание при этом, что матрица, находящаяся в произведении слева, соответствует матрице обратного теоретико-числового преобразования, т.е. выражению (3.2), а матрица, находящаяся справа от знака произведения, представляет прямое преобразование (3.1).

В (3.1), (3.2) было декларировано, что $N$ - простое число. Это связано с тем, что при простом $N$ в кольце вычетов $\mathcal{K}_m$ всегда существует обратный элемент $N^{-1}$, так как в этом случае образуется тройка взаимно простых чисел $s, N, M$ ($N$ -простое по определению, взаимная простота $s$ и $M$ следует из следствия 1 теоремы (3.1)). Если же $N$ не является простым числом, то возникают проблемы не только с обратным элементом $N^{-1}$, но и со значением выражения (3.4), для которого должны обязательно выполняться условия (3.7), определяющие существование преобразований (3.1), (3.2). С другой стороны произвольный выбор $N$ весьма желателен, так как он, в конечном счете, определяет размерность преобразования.

Пусть $N$ - составное число $N = p_1 p_2$, составленное из двух простых чисел $p_1$ и $p_2$. Тогда система весовых функций преобразования будет содержать степенные последовательности с периодами $N, p_1$ и $p_2$. Нетрудно убедиться в том, что средние значения этих функций при значениях $i = \overline{0, N-1}$ будут равны или кратны соответственно

$$S_1 = M = \sum_{m=0}^{N-1} s^m \text{ - период равен } N,$$

$$S_2 = M = \sum_{m=0}^{p_2-1} s^{p_1 m} \text{ - период равен } p_2,$$

$$S_3 = M = \sum_{m=0}^{p_1-1} s^{p_2 m} \text{ - период равен } p_1.$$

Если между значениями $p_1$ и $p_2$ выполняется соотношение $p_1 < p_2$, тогда между этими суммами существует соотношение

$$S_3 < S_2 < S_1. \qquad (3.8)$$

Для выполнения условия существования преобразований (3.7) необходимо, чтобы

$$S_1 = 0 \,|\, \mathrm{mod} M, \; S_2 = 0 \,|\, \mathrm{mod} M, \; S_3 = 0 \,|\, \mathrm{mod} M, \; \forall S_i > N.$$

Из (3.8) следует, что $M = S_3$, так как $S_3$ наименьшая сумма. Однако, теперь необходимо, чтобы $S_1$ и $S_2$ делились на величину $S_3$ без остатка. Действительно, на основании свойств геометрической прогрессии имеем

$$S_1 = \frac{s^N - 1}{s - 1}, \quad S_2 = \frac{(s^{p_1})^{p_2} - 1}{s^{p_1} - 1} = \frac{s^N - 1}{s^{p_1} - 1}, \quad S_3 = \frac{(s^{p_2})^{p_1} - 1}{s^{p_2} - 1} = \frac{s^N - 1}{s^{p_2} - 1}.$$

Тогда
$$\frac{S_1}{S_3} = \frac{s^N - 1}{s - 1} = \sum_{i=0}^{p_2 - 1} s^i, \qquad (3.9)$$

$$\frac{S_2}{S_3} = \frac{s^{p_2} - 1}{s^{p_1} - 1} = \sum_{i=0}^{p_2 - 1} s^i. \qquad (3.10)$$

говорит о том, что $N$ не может быть составным числом, так как не будут удовлетворяться условия (3.7). Исключение составляет случай, вытекающий из (3.9), при котором

$$N = p \cdot p. \qquad (3.11)$$

По индукции можно показать, что $N = p^n$. Тогда условие, связующее $p, s,$ и $M$ будет следующее:

$$M = \sum_{m=0}^{p-1} s^{pm} > p^n.$$

На основании свойств геометрической прогрессии

$$\frac{s^p - 1}{s - 1} > p^n \quad \text{или} \quad p < \sqrt[n]{\frac{s^p - 1}{s - 1}}. \qquad (3.12)$$

Таким образом, получено еще одно выражение для размерности преобразования $N$, при котором существует теоретико-числовое преобразование. Далее уточним величину модуля $M$. На основании (3.11) имеем

$$N = p^n = p \cdot p^{n-1}.$$

Тогда
$$\mathrm{I} \quad s - 1 = (s^{p^{n-\ldots 1}})^p - 1, \qquad (3.13)$$

или в силу (3.11)
$$s^N - 1 = (s^{p^{n-1}} - 1) \cdot \sum_{m=0}^{p-1} (s^{p^{n-1}})^m.$$

Однако
$$(s^{p^{n-1}} - 1) = (s - 1) \cdot \sum_{m=0}^{p^{n-1} - 1} s^k.$$

Подставляя в (3.13), получим
$$s^N - 1 = (s - 1) \cdot \sum_{k=0}^{p^{n-1} - 1} s^k \cdot \sum_{m=0}^{p-1} (s^{p^{n-1}})^m.$$

Поделив обе части полученного выражения на $(s - 1)$, имеем

$$M = \sum_{m=0}^{n-1} s^m = \sum_{k=0}^{p^{n-1} - 1} s^k \cdot \sum_{m=0}^{p-1} (s^{p^{n-1}})^m \qquad (3.14)$$

Таким образом, при размерности преобразования $N = p^n$ модуль $M$ является составным числом и в качестве модуля преобразования может быть выбрано одно из этих составных чисел. Вместе с тем условия преобразования удовлетворит только модуль $M$, равный

$$M = \sum_{m=0}^{p-1} (s^{p^{n-1}})^m \quad . \tag{3.15}$$

Этот факт объясняется тем, что в силу теоремы (3.1) модуль преобразования

$$M = \sum_{k=0}^{p^{n-1}-1} s^k \quad , \tag{3.16}$$

дает последовательность $s^k$ размерностью $p^{n-1}$, следовательно, такую же размерность преобразования, что значительно меньше требуемого $N$. Далее, на основании (3.15) ясно, что последовательность $s^k$ при изменении $i = \overline{0, N-1}$ и

$$s^i < M = \sum_{m=0}^{p-1} (s^{p^{n-1}})^m \tag{3.17}$$

имеют обычную явную степенную зависимость.

Определим теперь, какую зависимость будет представлять $s^i$, если $s^i > M$. Действительно, последний член в (3.15) имеет вид

$$(s^{p^{n-1}})^{p-1} = s^{p^n - p^{n-1}}.$$

Стало быть, элементы последовательности $s^i > M$ будут иметь выражение вида

$$s^i = (s^{p^n - p^{n-1}}) \cdot s^k = s^{p^n - p^{n-1} + k},$$

где, $k = \overline{1, N - p^n + p^{n-1} - 1}$ получено заменой переменной

$$i = p^n - p^{n-1} + k > p^n - p^{n-1} \quad . \tag{3.18}$$

Тогда $\quad s^{p^n - p^{n-1} + k} = s^{p^n} \cdot s^{-p^{n-1} + k}$

на основании теоремы (3.1) дает $s^{p^n - p^{n-1} + k} \xlongequal{\mathscr{X}_m} s^{-p^{n-1} + k}$.

Заменив обратно $k$ на $i$, из (3.18) получим $\quad s^{p^n - p^{n-1} + k} \xlongequal{\mathscr{X}_m} s^{-p^{n-1} - p^n + p^{n-1} + i} \xlongequal{\mathscr{X}_m} s^{-p^n + i}$ или окончательно $\quad s^i \xlongequal{\mathscr{X}_m} (s^{p^n - i})^{-1} \xlongequal{\mathscr{X}_m} (s^{N-1})^{-1}$, при $s^i > M$. \qquad (3.19)

Таким образом, элементы последовательности $s^i$ при $s^i > M$ представляют собой обратные элементы по отношению к элементам $s^i > M$. Обратимость при этом существует только в кольце $\mathcal{K}_m$. В результате **S**-преобразование (3.1), (3.2) может рассматриваться, с учетом индексов по (3.18), как косо симметричное двухстороннее преобразование вида

$$X(k) \stackrel{\mathcal{K}_m}{=\!=\!=} \sum_{i=-(p^{n-1}-1)}^{p^n-p^{n-1}} x(i) s^{-(ki) \bmod [-p^n+1, p^n-p^{n-i}]}, \qquad (3.20)$$

$$x(i) \stackrel{\mathcal{K}_m}{=\!=\!=} \sum_{k=(p^{n-1}-1)}^{p^n-p^{n-1}} X(k) s^{(ki) \bmod [-p^n+1, p^n-p^{n-i}]}, \qquad (3.21)$$

$$M = \sum_{m=0}^{p-1} (s^{p^{n-1}})^m. \qquad (3.22)$$

При этом периодичность весовой функции преобразования определяется не величиной $N$ (замкнутым интервалом $[0, N-1]$), а замкнутым интервалом $[-p^{n-1}+1, p^n-p^{n-1}]$, который к тому же и несимметричен, так как абсолютные значения его концов не равны между собой

$$|-p^{n-1}+1| \neq |p^n - p^{n-1}|,$$

отчего и название преобразования – косо симметричное.

Вместе с тем существует еще один уникальный случай, представляющий исключение из условия (3.11). Речь идет о четном значении $N$

$$N = 2p, \text{ где } p \text{ - любое целое число} \qquad (3.23)$$

Стало быть $\qquad s^N - 1 = s^{2p} - 1 = (s^p)^2 - 1 = (s^p - 1)(s^p + 1)$

Далее, в силу (3.11) $\qquad s^p - 1 = (s-1) \cdot \sum_{n=0}^{p-1} s^n.$

Тогда $\qquad s^N - 1 = (s-1) \cdot \left( \sum_{n=0}^{p-1} s^n \right) \cdot (s^p + 1).$

Разделив обе части на величину $(s-1)$, имеем

$$M = \sum_{m=0}^{p-1} (s^{p^{n-1}})^m = (s^p + 1) \cdot \sum_{n=0}^{p-1} s^n. \qquad (3.24)$$

Как и в случае (3.14), преобразованию удовлетворяет модуль

$$M = s^p + 1, \qquad (3.25)$$

дающий преобразование размерности $N$. При этом **S**-преобразование принимает вид

$$X(k) \stackrel{\mathcal{K}_m}{=\!=\!=} \sum_{i=-(p-1)}^{p} x(i) s^{-(ki) \bmod [-p+1, p]}, \qquad (3.26)$$

$$x(i) \stackrel{\mathcal{K}_m}{=\!=\!=} \frac{1}{N} \sum_{k=-(p-1)}^{p} X(k) s^{(ki) \bmod [-p+1, p]}, \qquad (3.27)$$

$$M = s^p + 1. \qquad (3.28)$$

с интервалом существования $\qquad [-p+1, p], \qquad (3.29)$

представляющим собой также разновидность теперь уже почти симметричного двухстороннего преобразования, поскольку значения весовой функции $s^{ki}$ при $s^{ki} > M$, как и в случае (3.19), являются обратными элементами в кольце $\mathcal{K}_m$ к значениям весовой функции $s^{ki}$ при $s^{ki} > M$.

Однако может так оказаться, что $s^p + 1$ - также составное число, что следует из теоремы Ферма:

$$s^{Np} - 1 = 0 \mid \mod(N+1). \tag{3.30}$$

Вместе с тем, все простые числа, за исключением числа 2, являются нечетными числами, а поскольку степенная функция не содержит членов, равных нулю, то ограниченная сверху модулем, равным $M = N+1$, будет содержать на своем главном периоде всего N различимых между собой значений, где $N$ - четное число.

Таким образом, мы имеем дело с преобразованием, размерность которого представлена четным числом. Очевидно, что речь идет о преобразовании (3.26), (3.27) и (3.28), модуль которого может быть составным числом, содержим в себе число $N+1$. Действительно, в силу (3.30), число $s^N - 1$ делится на число $N+1$ без остатка. Это означает, что один из сомножителей из (3.24) обязательно делится на $N+1$. Если же для (3.28) имеет место вышеизложенное, т.е. $M = s^p + 1$ - составное число, в которое входит $N+1$, то существует следующее преобразование:

$$X(k) \stackrel{\mathcal{K}_m}{=\!=} \sum_{i=0}^{N-1} x(i) s^{-(ki) \mod N}, \tag{3.31}$$

$$x(i) \stackrel{\mathcal{K}_m}{=\!=} \frac{1}{N} \sum_{k=0}^{N-1} X(k) s^{(ki) \mod N}, \tag{3.32}$$

$$M = N + 1. \tag{3.33}$$

При этом, несмотря на малую величину модуля $M$, всегда существует обратный элемент для числа $N$

$$N^2 - 1 = (N+1)(N-1).$$

Это дает сравнение по модулю $N+1$ следующего вида:

$$N^2 - 1 \equiv 0 \mid \mod(N+1).$$

Или $$N^2 \equiv 1 \mid \mod(N+1)$$

Отсюда надо полагать, что $$N \cdot N^{-1} \stackrel{\mathcal{K}_m}{=\!=} 1,$$

Тогда $$N \cdot N^{-1} \stackrel{\mathcal{K}_m}{=\!=} N^2 \stackrel{\mathcal{K}_m}{=\!=} N \cdot N.$$

Стало быть, обратным элементом к числу $N$ является само число $N$.

Весовая функция преобразования в этом случае будет содержать $N$ различных значений, которые лежат в замкнутом интервале $[1, N]$, в силу модуля (3.33), и представляют собой целые числа от 1 до $N$.

Тогда матрицу как прямого, так и обратного преобразования можно упорядочить путем перестановки строк и столбцов таким образом, что элементы строки с индексом $k=1$ ($i=1$) или столбца с индексом $i=1$ ($k=1$), соответственно матриц прямого и обратного преобразований, будут представлять линейно возрастающую последовательность, в которой разность между соседними элементами постоянна и равна 1. В этом случае будет иметь место одна из разновидностей пилообразного преобразования (весовая функция, при $k=1$ - для прямого и при $i=1$ - для обратного преобразований, представляет собой зависимость в виде отрезка прямой $y=x$, именуемой в технике пилой).

В заключение заметим, что часто используемые теоретико-числовые преобразования, такие как Мерсена и Ферма, являются частными случаями **S**-преобразования при $s=2$. Так, преобразование Мерсена получается из (3.1) при $s=2$ и модуле преобразования

$$M = 2^N - 1 = (2-1) \cdot \sum_{m=0}^{N-1} 2^m = \sum_{m=0}^{N-1} 2^m.$$

Преобразование же (3.26), (3.27), (3.28) дает преобразование Ферма при следующих условиях: $s=2$ и $p=2^n$, где $n$ - целое число.

### 4. Свойства весовых функций **S**-преобразования.

Заметим, что здесь и далее для кольца вычетов в $\mathcal{X}_m$ модуль может быть следующим

$$M = \sum_{m=0}^{N-1} s^m, \quad N\text{-простое число},$$

$$M = \sum_{m=0}^{p-1} s^{pm}, \quad N = p^n\text{-составное число}. \tag{4.1}$$

4.1. Функции **S**-преобразования периодичны. Их период равен $N$.

$$S(i, k+N) \stackrel{\mathcal{X}_m}{=\!=\!=} S(ik). \tag{4.2}$$

Действительно $\quad s^{i(k+N)\bmod N} = s^{(ik+iN)\bmod N} = s^{(ik)\bmod N}$

или в силу (4.1) имеем (4.2). Периодичность функций **S**-преобразования для двух значений переменной $i$ ($i=1$, $i=2$) иллюстрируется ниже поданным рисунком.

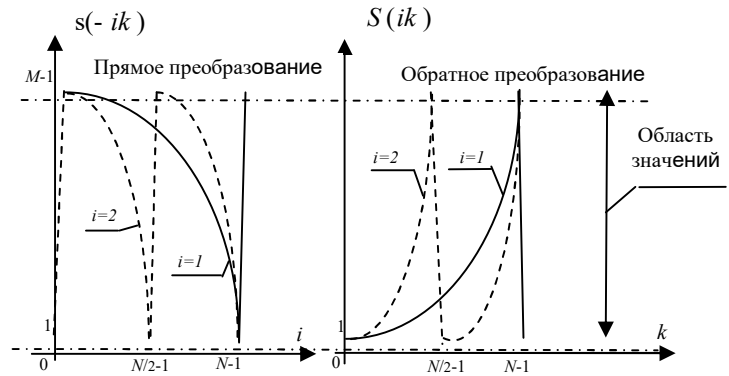

Рис. 4.1. График функций $S(1k)$ и $S(2k)$ $S$-преобразования.

4.2. Функции **S**-преобразования мультипликативные в кольце $\mathcal{K}_m$.

$$|S(a,k)S(b,k)| \stackrel{\mathcal{K}_m}{=\!=} S|(a+b)k|. \qquad (4.3)$$

Действительно $\quad S|a,k|\mathrm{mod}N \cdot S|b,k|\mathrm{mod}N = S|(a+b),k|\mathrm{mod}N$, в силу (4.1) имеем (4.3).

4.3. Среднее значение функций $S(ik)$ в $\mathcal{K}_m$ при $i \neq 0$ равно нулю.

$$\sum_{k=0}^{N-1} s(i,k) \stackrel{\mathcal{K}_m}{=\!=} 0, \quad \sum_{i=0}^{N-1} S(i,k) \stackrel{\mathcal{K}_m}{=\!=} 0, \qquad (4.4)$$

что следует из (3.2.7).

4.4. Система функций $S(ik)$ ортогональна в кольце $\mathcal{K}_m$

$$\sum_{k=0}^{N-1} S(a,k) \cdot S(b,k) \stackrel{\mathcal{K}_m}{=\!=} \sum_{k=0}^{N-1} [S(a+b,k)] \stackrel{\mathcal{K}_m}{=\!=} 0 \qquad (4.5)$$

Действительно

$$\sum_{k=0}^{N-1} s^{[ak]\mathrm{mod}\,N} \cdot s^{[bk]\mathrm{mod}\,N} \stackrel{\mathcal{K}_m}{=\!=} \sum_{k=0}^{N-1} s^{[(a+b)k]\mathrm{mod}\,N} \stackrel{\mathcal{K}_m}{=\!=} s^{[(a+b)k]\mathrm{mod}\,N} \sum_{k=0}^{N-1} s^{[k]\mathrm{mod}\,N} \stackrel{\mathcal{K}_m}{=\!=} 0.$$

В силу (4.1) имеем (4.5).

4.5. Матрицы преобразования, составленные из этих функций, симметричны относительно главной диагонали: $\quad |S(i,k)| \stackrel{\mathcal{K}_m}{=\!=} |S(k-i)|, \qquad (4.6)$

На основании перечисленных свойств функций можно утверждать, что **S**-преобразование в кольце $\mathcal{K}_m$ является ортогональным преобразованием с мультипликативным базисом. Матрицы преобразования могут быть факторизованы, следовательно, возможно построение быстрых алгоритмов расчета коэффициентов преобразования.

## 5. Основные теоремы S-преобразования.

При изложении теорем мы по-прежнему будем пользоваться терминологическими соглашениями об оригинале и изображении, чтобы не путать их со спектрами Фурье.

5.1. *Теорема о сумме (разности).*

Поэлементная (прямая) сумма (разность) последовательностей оригиналов в кольце вычетов приводит к сумме (разности) их изображений:

$$x(i) + y(i) \xrightarrow{\mathcal{X}_m} X(k) + Y(k) \qquad (5.1)$$

3.4.2. *Теорема линейности.*

Умножение последовательности-оригинала на постоянный множитель в кольце вычетов $\mathcal{X}_m$ приводит к умножению последовательности-изображения на тот же множитель:

$$\left[\sum_{m=o}^{N-1} \mu x(i) \cdot s^{ik}\right] \bmod M = \left[\mu \sum_{i=o}^{N-1} x(i) \cdot s^{ik}\right] \bmod M = \mu X(k), \qquad (5.2)$$

$$\mu x(i) \xrightarrow{\mathcal{X}_m} \mu x(k) \quad \text{при } \mu \in S.$$

5.3. *Теорема о сдвиге (теорема смещения).*

Поэлементное (прямое) перемножение последовательности-оригинала на весовую функцию **S**-преобразования в кольце вычетов $\mathcal{X}_m$ приводит к сдвигу последовательности-изображения:

$$X_C(k) \xrightarrow{\mathcal{X}_m} \sum_{i=0}^{N-1} x(i) s^{\pm ci} \cdot s^{-[ik] \bmod N} \xrightarrow{\mathcal{X}_m} \sum_{i=0}^{N-1} x(i) \cdot s^{-[(\pm c+k)i] \bmod N},$$

$$x(i) s^{-[(\pm c+k)i] \bmod N} \xrightarrow{\mathcal{X}_m} X[(\pm c+k) \bmod N]. \qquad (5.3)$$

Из выражения (5.3) видно, что сдвиг изображения является циклическим.

5.4. *Теорема о свертке изображения.*

Поэлементное произведение двух последовательностей-оригиналов в кольце вычетов $\mathcal{X}_m$ приводит к свертке их изображений:

$$x(i) \cdot y(i) \xrightarrow{\mathcal{X}_m} X(k) * Y(k), \qquad (5.4)$$

где * - (здесь и ниже) обозначение свертки, как оператора, а знак $\xrightarrow{\mathcal{X}_m}$ означает: влечет за собой, приводит к. Действительно, изображение от произведения $x(i) \cdot y(i)$

$$X(n) \xrightarrow{\mathcal{X}_m} \sum_{i=0}^{N-1} x(i) \cdot y(i) \cdot s^{-[in] \bmod N}.$$

Однако $\qquad x(i) \xrightarrow{\mathcal{X}_m} \dfrac{1}{N} \sum_{k=0}^{N-1} X(k) \cdot s^{-[ik] \bmod N}.$

Тогда $\qquad X(n) \xrightarrow{\mathcal{X}_m} \sum_{m=0}^{N-1} \dfrac{1}{N} \sum_{i=0}^{N-1} \left(X(k) \cdot s^{-[ik] \bmod N}\right) \cdot y(i) s^{-[in] \bmod N}.$

Изменив порядок суммирования, получим:

$$X(n) \xrightarrow{\mathcal{S}_m} \frac{1}{N}\sum_{k=0}^{N-1} X(k) \sum_{i=0}^{N-1} y(i) \cdot s^{-[(n-k)i]\mathrm{mod}N} \xrightarrow{\mathcal{S}_m} \frac{1}{N}\sum_{k=0}^{N-1} X(k) \cdot Y[(n-k)]\mathrm{mod}N.$$

Поскольку индекс $n-k$ находится в пределах модуля $N$, то полученное выражение описывает циклическую (круговую) свертку и таким образом

$$x(i) \cdot y(i) \xrightarrow{\mathcal{S}_m} \frac{1}{N}\sum_{k=0}^{N-1} X(k) \cdot Y[(n-k)\mathrm{mod}N]. \qquad (5.5)$$

Перейдем теперь к обратным утверждениям этих теорем.

5.5. *Поэлементная (прямая) сумма (разность)* последовательностей изображений приводит к сумме (разности) оригиналов:

$$X(k) + Y(k) \xrightarrow{\mathcal{S}_m} x(i) + y(i). \qquad (5.6)$$

5.6. *Умножение изображения* на постоянный множитель в кольце $\mathcal{S}_m$ приводит к умножению последовательности-оригинала на тот же множитель:

$$\mu \cdot X(k) \xrightarrow{\mathcal{S}_m} \mu \cdot X(l),\ \mu \in \mathbf{S}. \qquad (5.7)$$

5.7. *Циклический сдвиг последовательности-изображения* приводит к умножению в кольце $\mathcal{S}_m$ последовательности-оригинала на весовую функцию преобразования:

$$X[(\pm n + k)\mathrm{mod}N] \xrightarrow{\mathcal{S}_m} x(i) \cdot s^{\pm(n)\mathrm{mod}N}. \qquad (5.8)$$

5.8. *Циклическая свертка последовательностей-изображений* в кольце $\mathcal{S}_m$ влечет за собой поэлементное произведение последовательностей-оригиналов в том же кольце:

$$X(k) * Y(k) \xrightarrow{\mathcal{S}_m} \sum_{k=0}^{N-1} X(n) \cdot Y[(n-k)\mathrm{mod}N] \xrightarrow{\mathcal{S}_m} x(i) \cdot y(i). \qquad (5.9)$$

Покажем теперь, что имеют место и обратные теоремы S-преобразования.

5.9. Умножение в кольце $\mathcal{S}_m$ последовательности-изображения на весовую функцию преобразования приводит к циклическому сдвигу последовательности-оригинала:

$$X(k) \cdot s^{\pm[nk]\mathrm{mod}N} \xrightarrow{\mathcal{S}_m} x[(i \pm n)\mathrm{mod}N]. \qquad (5.10)$$

Действительно

$$x_c(i) \xrightarrow{\mathcal{S}_m} \frac{1}{N}\sum_{k=0}^{N-1} X(k) \cdot s^{\pm[nk]\mathrm{mod}N} \cdot s^{[ki]\mathrm{mod}N} \xrightarrow{\mathcal{S}_m} \frac{1}{N}\sum_{k=0}^{N-1} X(k) \cdot s^{[\pm(n-i)k]\mathrm{mod}N} \xrightarrow{\mathcal{S}_m} x[(i \pm n)\mathrm{mod}N],$$

$$X(k) \cdot s^{[\pm(n-i)k]\mathrm{mod}N} \xrightarrow{\mathcal{S}_m} x[(i \pm n)\mathrm{mod}N].$$

5.10. *Теорема о свертке оригинала.*

Прямое (поэлементное) произведение двух последовательностей-изображений в кольце $\mathcal{S}_m$ приводит к циклической свертке их оригиналов в том же кольце:

$$X(k) \cdot Y(k) \xrightarrow{\mathcal{S}_m} x(i) * y(i). \qquad (5.11)$$

Действительно, оригинал от произведения изображений

$$x(n) \xlongequal{\mathcal{K}_m} \frac{1}{N} \sum_{k=0}^{N-1} X(k) \cdot Y(k) \cdot s^{[nk] \bmod N},$$

где $x(n)$ -последовательность-изображение, представляющая свертку $x(i)$ и $y(i)$.

Однако 
$$X(k) \xlongequal{\mathcal{K}_m} \sum_{i=0}^{N-1} x(i) \cdot s^{-[ki] \bmod N}$$

$$x(n) \xlongequal{\mathcal{K}_m} \frac{1}{N} \sum_{k=0}^{N-1} Y(k) \left( \sum_{i=0}^{N-1} x(i) \cdot s^{-[ki] \bmod N} \right) s^{[nk] \bmod N},$$

Изменив порядок суммирования, имеем

$$x(n) \xlongequal{\mathcal{K}_m} \sum_{k=0}^{N-1} x(i) \left( \frac{1}{N} \sum_{i=0}^{N-1} Y(k) \cdot s^{[(n-i)k] \bmod N} \right) \xlongequal{\mathcal{K}_m} \sum_{i=0}^{N-1} x(i) \cdot y[(n-i) \bmod N].$$

Последнее и есть циклическая свертка

$$x(i) * y(i) \xlongequal{\mathcal{K}_m} \sum_{i=0}^{N-1} x(i) \cdot y[(n-i) \bmod N]. \tag{5.12}$$

Отсюда следует, что имеют место и обратные теоремы.

5.11. *Теорема Винера-Хинчина.*

Автокорреляционная функция оригинала и плотность энергии изображения в кольце $\mathcal{K}_m$ связаны между собой S-преобразованием

$$\sum_{i=0}^{N-1} x(i) \cdot y[(n-i) \bmod N] \xlongequal{\mathcal{K}_m} \frac{1}{N} \sum_{i=0}^{N-1} Y(k)^2 \cdot s^{(kn) \bmod N}. \tag{5.13}$$

*Доказательство*:

Полагая в (5.11) $x(i) = y(i)$, $X(k) = Y(k)$, имеем

$$Y(k) \cdot Y(k) \xrightarrow{\mathcal{K}_m} \sum_{i=0}^{N-1} y(i) \cdot y[(n-i) \bmod N]$$

или 
$$Y(k)^2 \xrightarrow{\mathcal{K}_m} \sum_{i=0}^{N-1} y(i) \cdot y[(n-i) \bmod N] \tag{5.14}$$

Таким образом, плотность энергии изображения представляет собой свертку оригинала самого с собой или автокорреляцию (автоковариацию) в кольце $\mathcal{K}_m$. Установим связь между плотностью энергии изображения и автокорреляцией оригинала. Из (3.2) следует

$$y(i) \xlongequal{\mathcal{K}_m} \frac{1}{N} \sum_{k=0}^{N-1} Y(k) \cdot s^{[ik] \bmod N}.$$

Заменив $i$ на $n-i$, получим 
$$y(n-i) \xlongequal{\mathcal{K}_m} \frac{1}{N} \sum_{k=0}^{N-1} Y(k) \cdot s^{[(n-i)k] \bmod N}.$$

Тогда $$\sum_{i=0}^{N-1} y(i) \cdot y[(n-i) \bmod N] \stackrel{Z_m}{=\!=\!=} \sum_{i=0}^{N-1} y(i) \left( \frac{1}{N} \sum_{k=0}^{N-1} Y(k) \cdot s^{[(n-i)k] \bmod N} \right).$$

Или $$\sum_{i=0}^{N-1} y(i) \cdot y[(n-i) \bmod N] \stackrel{\mathscr{X}_m}{=\!=\!=} \sum_{i=0}^{N-1} y(i) \left( \frac{1}{N} \sum_{k=0}^{N-1} Y(k) \cdot s^{[nk] \bmod N} \cdot s^{-[ik] \bmod N} \right).$$

Изменив порядок суммирования, получим:

$$\sum_{i=0}^{N-1} y(i) \cdot y[(n-i) \bmod N] \stackrel{\mathscr{X}_m}{=\!=\!=} \frac{1}{N} \sum_{k=0}^{N-1} Y(k) \cdot s^{[nk] \bmod N} \sum_{i=0}^{N-1} y(i) \cdot s^{-[ik] \bmod N}.$$

Внутренняя сумма представляет изображение от последовательности $y(i)$, и равна $Y(k)$, и

тогда $$\sum_{i=0}^{N-1} y(i) \cdot y[(n-i) \bmod N] \stackrel{\mathscr{X}_m}{=\!=\!=} \frac{1}{N} \sum_{i=0}^{N-1} Y(k) \cdot Y(k) s^{(kn) \bmod N}.$$

или $$\sum_{i=0}^{N-1} y(i) \cdot y[(n-i) \bmod N] \stackrel{\mathscr{X}_m}{=\!=\!=} \frac{1}{N} \sum_{k=0}^{N-1} Y^2(k) \cdot s^{(kn) \bmod N}.$$

Итак, получено выражение (3.4.13), что и требовалось доказать. Можно сформулировать и обратное утверждение этой теоремы.

5.12. Автокорреляционная функция изображения связана с энергией оригинала в кольце $\mathscr{X}_m$ S-преобразованием

$$\sum_{k=0}^{N-1} Y(k) \cdot Y[(n-k) \bmod N] \stackrel{\mathscr{X}_m}{=\!=\!=} \frac{1}{N} \sum_{i=0}^{N-1} y^2(i) s^{-(ki) \bmod N}. \qquad (5.15)$$

Полагая, что в (3.3.4) $Y(k) \stackrel{\mathscr{X}_m}{=\!=\!=} X(k), y(i) \stackrel{\mathscr{X}_m}{=\!=\!=} x(i)$, получим

$$y(i) \cdot y(i) \xrightarrow{Z_m} \frac{1}{N} \sum_{k=0}^{N-1} Y(k) \cdot Y[(n-k) \bmod N],$$

или $$y(i)^2 \xrightarrow{\mathscr{X}_m} \frac{1}{N} \sum_{k=0}^{N-1} Y(k) \cdot Y[(n-k) \bmod N],$$

что далее из (3.2.1) следует $$Y(k) \stackrel{\mathscr{X}_m}{=\!=\!=} \sum_{i=0}^{N-1} y(i) \cdot s^{-(ki) \bmod N}.$$

Заменив $k$ на $n-i$, имеем $$Y[(n-k) \bmod N] \stackrel{\mathscr{X}_m}{=\!=\!=} \sum_{i=0}^{N-1} y(i) \cdot s^{-[(n-k)i] \bmod N}.$$

Тогда $$\frac{1}{N} \sum_{k=0}^{N-1} Y(k) \cdot Y[(n-k) \bmod N] \stackrel{\mathscr{X}_m}{=\!=\!=} \frac{1}{N} \sum_{k=0}^{N-1} Y(k) \left( \sum_{i=0}^{N-1} y(i) s^{-[(n-k)i] \bmod N} \right).$$

Или $$\frac{1}{N} \sum_{k=0}^{N-1} Y(k) Y[(n-k) \bmod N] \stackrel{\mathscr{X}_m}{=\!=\!=} \frac{1}{N} \sum_{k=0}^{N-1} Y(k) \left( \sum_{i=0}^{N-1} y(i) s^{-[ni] \bmod N} s^{[ki] \bmod N} \right).$$

Изменив порядок суммирования, получим

$$\frac{1}{N} \sum_{k=0}^{N-1} Y(k) \cdot Y[(n-k) \bmod N] \stackrel{\mathscr{X}_m}{=\!=\!=} \frac{1}{N} \sum_{i=0}^{N-1} y(i) \cdot s^{-[ni] \bmod N} \sum_{k=0}^{N-1} Y(k) \cdot s^{[ki] \bmod N}.$$

Внутренняя сумма в правой части представляет собой $y(i)$, в силу чего имеем

$$\frac{1}{N} \sum_{k=0}^{N-1} Y(k) \cdot Y[(n-k) \bmod N] \stackrel{\mathscr{X}_m}{=\!=\!=} \frac{1}{N} \sum_{i=0}^{N-1} y(i) \cdot y(i) s^{-[ni] \bmod N},$$

или
$$\sum_{k=0}^{N-1} Y(k) \cdot Y[(n-k) \bmod N] \stackrel{\mathcal{K}_m}{=\!=\!=} \sum_{i=0}^{N-1} y(i)^2 \cdot s^{-[ni] \bmod N}.$$

Таким образом, плотность мощности оригинала связана с автокорреляцией изображения теоретико-числовым преобразованием, что является обратным утверждением аналога теоремы Винера - Хинчина.

5.13. *Равенство Парсеваля*.

Энергия оригинала равна энергии изображения:

$$\sum_{i=0}^{N-1} y(i)^2 \stackrel{\mathcal{K}_m}{=\!=\!=} \sum_{k=0}^{N-1} Y(k)^2, \qquad (5.16)$$

Применив обратное преобразование к двум сторонам выражения (3.4.15), получим

$$y(i)^2 \stackrel{\mathcal{K}_m}{=\!=\!=} \frac{1}{N} \sum_{n=0}^{N-1} \left( \sum_{k=0}^{N-1} Y(k) \cdot Y[(n-k) \bmod N] \right) s^{(ni) \bmod \ldots N},$$

а, суммируя по $i$ в обеих частях равенства, имеем

$$\sum_{i=0}^{N-1} y(i)^2 \stackrel{\mathcal{K}_m}{=\!=\!=} \frac{1}{N} \sum_{i=0}^{N-1} \left[ \sum_{n=0}^{N-1} \left( \sum_{k=0}^{N-1} Y(k) \cdot Y[(n-k) \bmod N] \right) s^{(ni) \bmod N} \right],$$

Изменим порядок суммирования:

$$\sum_{i=0}^{N-1} y(i)^2 \stackrel{\mathcal{K}_m}{=\!=\!=} \sum_{i=0}^{N-1} \sum_{k=0}^{N-1} Y(k) \left( \frac{1}{N} \sum_{n=0}^{N-1} Y(n-k) s^{(ni) \bmod N} \right),$$

Очевидно, что в силу теоремы о сдвиге внутренняя сумма

$$\frac{1}{N} \sum_{k=0}^{N-1} Y(n-k) \cdot s^{(ni)} \stackrel{\mathcal{K}_m}{=\!=\!=} y(i) \cdot s^{-[ki] \bmod N}.$$

Тогда
$$\sum_{i=0}^{N-1} y(i)^2 \stackrel{\mathcal{K}_m}{=\!=\!=} \sum_{i=0}^{N-1} \sum_{k=0}^{N-1} Y(k) \cdot y(i) \cdot s^{-[ki] \bmod N},$$

и, еще раз изменив порядок суммирования, получим

$$\sum_{i=0}^{N-1} y(i)^2 \stackrel{\mathcal{K}_m}{=\!=\!=} \sum_{k=0}^{N-1} Y(k) \sum_{i=0}^{N-1} y(i) \cdot s^{-[ki] \bmod N},$$

Снова внутренняя сумма равна $Y(k)$, тогда

$$\sum_{i=0}^{N-1} y(i)^2 \stackrel{\mathcal{K}_m}{=\!=\!=} \sum_{k=0}^{N-1} Y(k) \cdot Y(k),$$

или
$$\sum_{i=0}^{N-1} y(i)^2 \stackrel{\mathcal{K}_m}{=\!=\!=} \sum_{k=0}^{N-1} Y(k)^2,$$

что и требовалось доказать.

## 6. Связь между S-преобразованиями в различных базисах.

Под различными базисами S-преобразования будем понимать то, что в выражениях (3.2) и (3.1) значения $s$ могут принимать различные значения, т.е. $s_1 \neq s_2$. Тогда возникает необходимость пересчета результатов из одного базиса в другой.

Пусть последовательность $x(i)$ в базисе, основание которого составляет число $s_2$, имеет изображение $X(k)$ и требуется найти изображение $X(n)$ той же последовательности-оригинала в базисе $s_1$. Изображение последовательности $x(i)$ в базисе $s_1$

$$X(n) \stackrel{\mathscr{K}_m}{=\!=\!=} \sum_{i=0}^{N-1} x(i) \cdot s_1^{-(in) \bmod N},$$

однако
$$x(i) \stackrel{\mathscr{K}_m}{=\!=\!=} \frac{1}{N} \sum_{k=0}^{N-1} X(k) \cdot s_2^{(ki) \bmod N}.$$

Тогда
$$X(n) \stackrel{\mathscr{K}_m}{=\!=\!=} \sum_{i=0}^{N-1} \left( \frac{1}{N} \sum_{k=0}^{N-1} X(k) \cdot s_2^{[ki] \bmod N} \right) s_1^{-(in) \bmod N}.$$

Изменив порядок суммирования, получим

$$X(n) \stackrel{\mathscr{K}_m}{=\!=\!=} \frac{1}{N} \sum_{k=0}^{N-1} X(k) \sum_{i=0}^{N-1} s_1^{-(in) \bmod N} \cdot s_2^{[ki] \bmod N}.$$

Внутренняя сумма является функцией двух переменных $n$ и $k$, представляющей собой не что иное, как ядро S-преобразования:

$$\Phi(k,n) \stackrel{\mathscr{K}_m}{=\!=\!=} \sum_{i=0}^{N-1} s_2^{[ki] \bmod N} \cdot s_1^{-(in) \bmod N}$$

на основании (3.1) $\Phi(k,n) \stackrel{\mathscr{K}_m}{=\!=\!=} \sum_{i=0}^{N-1} [(s_2^k)^i \bmod M_2] \cdot [(s_1^{-n})^i \bmod M_1]$

$$\text{где} \quad M_1 = \sum_{m=0}^{N-1} s_1^m, M_2 = \sum_{m=0}^{N-1} s_2^m.$$

Далее, $M = M_1 \cdot M_2$ и $\quad \Phi(k,n) \stackrel{\mathscr{K}_m}{=\!=\!=} (\sum_{i=0}^{N-1} (s_2^k \cdot s_1^{-n})^i) \bmod M$. (6.1)

при этом $\quad M = \left[ \sum_{m=0}^{N-1} s_2^m \right] \cdot \left[ \sum_{m=0}^{N-1} s_1^m \right] = \sum_{m=0}^{N-1} s_2^m \sum_{m=0}^{N-1} s_1^m,$ (6.2)

$k \in [0, N-1], n \in [0, N-1]$. здесь $M$ - новый модуль кольца $\mathscr{K}_m$.

Окончательно $\quad X(n) \stackrel{\mathscr{K}_m}{=\!=\!=} \left[ \frac{1}{N} \sum_{k=0}^{N-1} X(k) \cdot \Phi(k,n) \right] \bmod M.$ (6.3)

Таким образом, изображение в кольце $\mathscr{K}_m$, где $M$ из (3.2), в базисах $s_1$ и $s_2$ связаны между собой ядром преобразования $\Phi(k,n)$. Аналогичным образом можно сформулировать обратное утверждение.

Пусть изображение $X(k)$ имеет оригинал $x(n)$ в базисе $s_1$. Требуется по оригиналу $x(n)$ найти оригинал $x(i)$ в базисе $s_2$. Оригинал в базисе $s_1$

$$x(i) \overset{\dot{\mathscr{X}}_m}{=\!=} \frac{1}{N}\sum_{k=0}^{N-1} X(k)\cdot s_1^{[ik]\bmod N}.$$

Однако
$$X(k) \overset{\dot{\mathscr{X}}_m}{=\!=} \sum_{n=0}^{N-1} x(n)\cdot s_2^{-(nk)\bmod N}.$$

Тогда
$$x(i) \overset{\dot{\mathscr{X}}_m}{=\!=} \frac{1}{N}\sum_{k=0}^{N-1}\bigl(\sum_{n=0}^{N-1} x(n)\cdot s_2^{-[nk]\bmod N}\bigr)\cdot s_1^{[ik]\bmod N}.$$

Изменив порядок суммирования, получим

$$x(i) \overset{\dot{\mathscr{X}}_m}{=\!=} \frac{1}{N}\sum_{n=0}^{N-1} x(n)\cdot \bigl(\sum_{k=0}^{N-1} s_2^{-(nk)\bmod N}\cdot s_1^{(ik)\bmod N}\bigr).$$

Внутренняя сумма является тоже ядром преобразования

$$\Phi(k,n) \overset{\dot{\mathscr{X}}_m}{=\!=} \sum_{i=0}^{N-1} s_2^{-(nk)\bmod N}\cdot s_1^{(ik)\bmod N}, \tag{6.4}$$

или на основании (3.1.1)
$$\Phi(k,n) \overset{\dot{\mathscr{X}}_m}{=\!=} \bigl(\sum_{i=0}^{N-1} (s_2^{-n}\cdot s_1^{i})^{k}\bigr). \tag{6.5}$$

где $M = \sum_{m=0}^{N-1} s_1^m \cdot \sum_{m=0}^{N-\dots 1} s_2^m$, - модуль кольца $\dot{\mathscr{X}}_m$.

Окончательно
$$x(i) \overset{\dot{\mathscr{X}}_m}{=\!=} \left[\frac{1}{N}\sum_{n=0}^{N-1} x(n)\cdot \Phi(n,i)\right]\bmod M_2.$$

Следовательно, два оригинала в разных базисах, имеющих одно и то же изображение, связаны в кольце $\dot{\mathscr{X}}_m$ между собой ядром преобразования $\Phi(n,i)$. В силу того, что матрица ядер **S**-преобразования несимметрична, очевидно, что $|\Phi(k,n)|\neq|\Phi(n,i)|$ даже при перестановке индексов.

### 7. Комплексное теоретико-числовое преобразование.

При определении теоретико-числовых преобразователей по выражениям (3.1) и (3.2) полагалось, что $s$ - некоторое число, в общем случае комплексное. Однако если $s$ - комплексное число, то эти преобразования имеют особенности, которые следует рассматривать подробнее.

Итак, пусть $\dot{s} = a + jb$ - комплексное число, в котором $a = \mathrm{Re}\,s$, $b = \mathrm{Im}\,s$ - взаимно простые числа. Покажем, что теорема (3.1) имеет место и для комплексных чисел.

$$\dot{s}^{(x)\bmod p} \overset{\dot{\mathscr{X}}_m}{=\!=} (\dot{s}^x)\bmod(\dot{s}^p - 1). \tag{7.1}$$

Первоначально по аналогии с (3.2) комплексный модуль кольца $\dot{\mathscr{X}}_m$ - $\dot{M} = \dot{s}^p - 1$, квадрат нормы которого

$$\|\dot M\| = \mathrm{Re}^2(\dot s^{\,p}-1) + \mathrm{Im}^2(\dot s^{\,p}-1). \qquad (7.2)$$

Тогда с одной стороны имеем $\quad [\dot c(\dot s^{\,p}-1) + (\dot s^{\,x})\,\mathrm{mod}\,(\dot s^{\,p}-1)]\,\mathrm{mod}\,(\dot s^{\,p}-1). \qquad (7.5)$

откуда следуют два сравнения
$$\begin{aligned}\dot s^{\,kp} &= 1\,|\,\mathrm{mod}(\dot s^{\,p}-1),\\ \dot c(\dot s^{\,p}-1) &= 0\,|\,\mathrm{mod}(\dot s^{\,p}-1).\end{aligned} \qquad (7.6)$$

Второе сравнение имеет место в силу линейности сравнений. При этом составлять систему двух вещественных сравнений не имеет смысла, так как из записи видно, что модуль сравнения равен одному из сомножителей в левой части этого сравнения. В первом же сравнении (7.6) к левой части прибавим и вычтем вещественную единицу:

$$\dot s^{\,kp} + 1 - 1 = 1\,|\,\mathrm{mod}\,(\dot s^{\,p}-1).$$

Далее, на основании теоремы Безу

$$\dot s^{\,kp} - 1 = \dot s^{\,p/k} - 1^k = (\dot s^{\,p}-1)\cdot F(\dot s^{\,p})^{k-1}.$$

Первое сравнение из (3.6.6) можно переписать таким образом:

$$(\dot s^{\,p}-1)\cdot F(\dot s^{\,p})^{k-1} + 1 \equiv 1\,|\,\mathrm{mod}\,(\dot s^{\,p}-1),$$

или $\quad \left[(\dot s^{\,p}-1)\cdot F(\dot s^{\,p})^{k-1}\right]\mathrm{mod}\,(\dot s^{\,p}-1) + (1)\mathrm{mod}\,(\dot s^{\,p}-1) \equiv 1\,|\,\mathrm{mod}\,(\dot s^{\,p}-1).$

Первое слагаемое в левой части этого сравнения сравнимо с нулем в силу того, что содержит сомножитель, равный модулю сравнения. Тогда $\quad (1)\mathrm{mod}(\dot s^{\,p}-1) \equiv 1\,|\,\mathrm{mod}\,(\dot s^{\,p}-1).$

Далее для корректности доказательства необходимо показать, что вещественное число, равное 1, имеет вычет по комплексному модулю $\dot s^{\,p}-1$, равный вещественной единице. Действительно, по определению вычетов по комплексному модулю [4] имеем два вещественных сравнения в кольце вычетов $\dot{\mathcal{X}}_m$

$$Re(\dot s^{\,p}-1)\,Re\,\dot x + Im(\dot s^{\,p}-1)\,Im\,\dot x \xlongequal{\dot{\mathcal{X}}_m} Re(\dot s^{\,p}-1)\cdot 1,$$

$$-Im(\dot s^{\,p}-1)\,Re\,\dot x + Re(\dot s^{\,p}-1)\,Im\,\dot x \xlongequal{\dot{\mathcal{X}}_m} Im(\dot s^{\,p}-1)\cdot 1,$$

где модуль кольца вычетов $\dot{\mathcal{X}}_m$ равен квадрату нормы комплексного числа $(\dot s^{\,p}-1)$, $\dot x = Re\,\dot x + j\,Im\,\dot x.$ - наименьший вычет.

Последнее утверждение можно представить в более удобной матричной форме

$$\begin{bmatrix} Re(\dot s^{\,p}-1) & Re(\dot s^{\,p}-1) \\ -Im(\dot s^{\,p}-1) & Re(\dot s^{\,p}-1) \end{bmatrix} * \begin{bmatrix} Re\,\dot x \\ Im\,\dot x \end{bmatrix} \xlongequal{\dot{\mathcal{X}}_m} \begin{bmatrix} Re(\dot s^{\,p}-1) \\ Im(\dot s^{\,p}-1) \end{bmatrix}$$

Решение этой системы в кольце вычетов $\dot{\mathcal{X}}_m$ дает $\dot x = 1 + j0$, следовательно, вычет равен вещественной единице. Стало быть, первое сравнение из (7.6) имеет место и тогда (7.5) можно переписать в виде

$$(1 \cdot \dot{S}^{(x) \bmod p}) \bmod (\dot{S}^p - 1) \stackrel{\dot{\mathcal{X}}_m}{=\!=} (0 + \dot{s}^x) \bmod (\dot{s}^p - 1),$$

или окончательно $\quad \dot{s}^{(x) \bmod p} \stackrel{\dot{\mathcal{X}}_m}{=\!=} (\dot{s}^x) \bmod (\dot{s}^p - 1).$

Таким образом, показано, что теорема (3.1) справедлива и для комплексной переменной $\dot{s}$, и теперь можно сформулировать комплексное теоретико-числовое преобразование.

Пусть задана некоторая комплексная функция (последовательность) в кольце вычетов $\dot{\mathcal{X}}_m$ вида $\dot{z}(i) = x(i) + jy(i)$, где $x(i)$ и $y(i)$ - вещественные последовательности, заданные в кольце вычетов $\dot{\mathcal{X}}_m$ и удовлетворяющие условиям $\begin{vmatrix} x(i) \end{vmatrix} < \|\dot{M}\|, \\ \begin{vmatrix} y(i) \end{vmatrix} < \|\dot{M}\|,$ а $N$ - простое вещественное число. Тогда пара комплексных преобразований примет вид

$$\dot{Z}(k) \stackrel{\dot{\mathcal{X}}_m}{=\!=} \sum_{i=0}^{N-1} \dot{z}(i) \dot{s}^{-(ki) \bmod N},$$

$$\dot{z}(i) \stackrel{\dot{\mathcal{X}}_m}{=\!=} \frac{1}{N} \sum_{k=0}^{N-1} \dot{Z}(k) \dot{s}^{(ki) \bmod N}, \quad (7.7)$$

$$\dot{M} = \sum_{m=0}^{N-1} \dot{s}^m.$$

Несмотря на совершенно новую форму записи, это преобразование при определенных фиксированных значениях $\dot{s}$ и $\dot{M}$ известно в литературе как преобразование Гаусса [5]. Это объясняется тем, что способ вычисления вычетов по комплексному модулю $\dot{M}$ осуществляется на основании первой и второй фундаментальных теорем Гаусса. Поэтому при непосредственном вычислении коэффициентов преобразования (7.7) необходимо решать в кольце вещественных вычетов $\dot{\mathcal{X}}_m$ $N$ - раз систему линейных алгебраических уравнений второго порядка вида

$$\begin{bmatrix} \operatorname{Re} \dot{M} & \operatorname{Im} \dot{M} \\ -\operatorname{Im} \dot{M} & \operatorname{Re} \dot{M} \end{bmatrix} * \begin{bmatrix} Re\, \dot{Z}(k) \\ Im\, \dot{Z}(k) \end{bmatrix} \stackrel{\dot{\mathcal{X}}_m}{=\!=} \begin{bmatrix} \operatorname{Re} \dot{M} & \operatorname{Im} \dot{M} \\ -\operatorname{Im} \dot{M} & \operatorname{Re} \dot{M} \end{bmatrix} * \begin{bmatrix} Re(\sum_{i=0}^{N-1} \dot{z}(i) \dot{s}^{-(ki) \bmod N}) \\ Im(\sum_{i=0}^{N-1} \dot{z}(i) \dot{s}^{-(ki) \bmod N}) \end{bmatrix}, \quad (7.8)$$

для прямого преобразования при $k = \overline{0, N-1}$ и, соответственно, для обратного преобразования при $i = \overline{0, N-1}$

$$\begin{bmatrix} \operatorname{Re} \dot{M} & \operatorname{Im} \dot{M} \\ -\operatorname{Im} \dot{M} & \operatorname{Re} \dot{M} \end{bmatrix} * \begin{bmatrix} Re\, \dot{Z}(k) \\ Im\, \dot{Z}(k) \end{bmatrix} \stackrel{\dot{\mathcal{X}}_m}{=\!=} \begin{bmatrix} \operatorname{Re} \dot{M} & \operatorname{Im} \dot{M} \\ -\operatorname{Im} \dot{M} & \operatorname{Re} \dot{M} \end{bmatrix} * \begin{bmatrix} Re(\frac{1}{N}\sum_{i=0}^{N-1} \dot{z}(i) \dot{s}^{-(ki) \bmod N}) \\ Im(\frac{1}{N}\sum_{i=0}^{N-1} \dot{z}(i) \dot{s}^{-(ki) \bmod N}) \end{bmatrix}. \quad (7.9)$$

При этом модуль кольца $\dot{\mathcal{X}}_m$ равен $\quad M = \|\dot{M}\| = \left\| \sum_{m=0}^{N-1} \dot{s}^m \right\|.$

## 8. Двойственность (дуализм) теорем **S**-преобразования

Под структурой или решеткой [6] (уже упоминалось во второй главе) понимают частично-упорядоченное множество, которое имеет точную верхнюю и точную нижнюю грань. Наиболее часто используют следующие соотношения:

$$\sup \mathbf{S} = \sup\{a,b\} = \max(a+b), \qquad (8.1)$$
$$\inf \mathbf{S} = \inf\{a,b\} = \min(a \cdot b),$$

при этом знак минимума может отсутствовать.

Однако такое определение точных граней множества не всегда удобно, так как связывает операции над элементами множества с выражениями (8.1). Для теоретико-числовых преобразований введем другое определение этих граней следующим образом:

$$\sup \mathbf{S} = \sup\{a,b\} = \max[(a*b) \bmod M] = M-1, \qquad (8.2)$$
$$\inf \mathbf{S} = \inf\{a,b\} = \min[\operatorname{int}_p(a*b)] = 0,$$

где $*$ - знак произвольной алгебраической операции, результат которой ограничен сверху операцией по модулю $M$,

$\operatorname{int}_p()$ - целая часть по отношению к числу $p$.

При $p=1$ - эта функция принимает традиционные представления, в другом случае, она отмечает тот факт, что целая часть существует не только при делении целых чисел, но и дробных.

Заметим, что функция $\operatorname{int}_p()$ может быть представлена также как

$$\operatorname{int}_p(A) = A - \frac{A}{p} \bmod M_1,$$

где $A$ - некоторое число или результат той или иной операции.

Из формулы (3.2) видно, что такое определение структуры очень хорошо согласуется с теорией вычетов по модулю $M$ и, следовательно, не имеет противоречий с вышеизложенной теорией **S**-преобразований. Далее, при изложении фундаментальной теоремы теоретико-числовых преобразований было использовано понятие кольца вычетов по модулю, как алгебры вида $\quad\mathcal{K}_m = <\mathbf{S, +, -, 0, 1}>.\qquad (8.3)$

При этом декларировалось, но никогда не использовалось соотношение $\mathbf{0 \neq 1}$, т.е. нейтральный элемент $\mathbf{0}$ операции сложения и нейтральный элемент $\mathbf{1}$ операции умножения не равны между собой, а, стало быть, различимы. Однако на практике это условие не всегда выполняется, например дополнительный код в ЭВМ и т.п. Если положить, что $\mathbf{0=1}$, т.е. нейтральные элементы операций сложения и умножения совпадают и равны единице, то возникает двойственность.

Действительно, пусть $a+a'=\mathbf{1}$, где $a'$ - обратный элемент для $a$ по отношению к операции сложения (дополнение до единицы). Тогда для структуры получим

$$a+a'=M+\mathbf{1}. \quad (8.4)$$

Далее, используя правила де - Моргана, имеем

$$(a'*b')' = M+\mathbf{1}-(M+\mathbf{1}-a)(M+\mathbf{1}-b).$$

Окончательно для кольца вычетов, заданного на структуре, получим двойственность вида

$$(a'*b')' = a+b - ab, \quad (8.5)$$

откуда следуют еще два выражения двойственности:

$$a*b = a+b-(a'*b')' \quad (8.6)$$

и

$$a+a = ab+(a'b')'. \quad (8.7)$$

Таким образом, установлено, что операции сложения и умножения, заданные на структуре по выражению (8.2), двойственны, т.е. могут быть определены двумя различными путями. Но так как операции заданные на структуре двойственны, то двойственны и теоремы S-преобразования. Так теорема о сумме (5.1) имеет следующие двойственные выражения, основанные на (8.7):

$$x(i)+y(i) \xrightarrow{\mathscr{Z}_m} X(k)*Y(k)+(X'(k)*Y'(k))', \quad (8.8)$$

$$x(i)*y(i)+(x'(i)*y'(i))' \xrightarrow{\mathscr{Z}_m} X(k)+Y(k), \quad (8.9)$$

$$x(i)*y(i)+(x'(i)*y'(i))' \xrightarrow{\mathscr{Z}_m} X(k)*Y(k)+(X'(k)*Y'(k))', \quad (8.10)$$

где $x'(i)$, $y'(i)$, $X'(k)$, $Y'(k)$ - дополнения до единицы в кольце вычетов $Z_m$ значений $x(i)$, $y(i)$, $X(k)$, $Y(k)$ соответственно.

Для теоремы линейности (5.2) тоже существуют двойственные ей выражения, получаемые из (3.7.6):

$$\mu'(x(i) \xrightarrow{\mathscr{Z}_m} \mu+X(k)-(\mu'*X'(k))', \quad (8.11)$$

$$\mu+x(i)-(\mu'*x'(i))' \xrightarrow{\mathscr{Z}_m} \mu X(k), \quad (8.12)$$

$$\mu+x(i)-(\mu'*x'(i))' \xrightarrow{\mathscr{Z}_m} \mu+X(k)-(\mu'*X'(k))'. \quad (8.13)$$

Интересную двойственность имеет теорема о свертке (3.5):

$$x(i)*y(i) \xrightarrow{\mathscr{Z}_m} \frac{1}{N}\sum_{k=0}^{N-1}X(k)+\frac{1}{N}\sum_{k=0}^{N-1}Y(n-k) - \frac{1}{N}\sum_{k=0}^{N-1}(X'(k)*Y'(n-k))', \quad (8.14)$$

$$x(i)+y(i)-(x'(i)*y'(i))' \xrightarrow{\mathscr{Z}_m} \frac{1}{N}\sum_{k=0}^{N-1}X(k)*Y(n-k), \quad (8.15)$$

$$x(i)+y(i)-(x'(i)*y'(i))' \xrightarrow{\mathscr{Z}_m} \frac{1}{N}(\sum_{k=0}^{N-1}X(k)+\sum_{k=0}^{N-1}Y(n-k)-\sum_{k=0}^{N-1}(X'(k)*Y'(n-k))'. \quad (8.16)$$

И, наконец, приведем двойственность равенства Парсеваля (5.16).

$$2\sum_{i=0}^{N-1}y(i) - \sum_{i=0}^{N-1}[(y'(i))^2]' \xrightarrow{Z_m} \sum_{k=0}^{N-1}Y^2(k), \quad (8.17)$$

$$\sum_{i=0}^{N-1} y^2(i) \xrightarrow{Z_m} 2\sum_{k=0}^{N-1} Y(k) - \sum_{i=0}^{N-1}[(y'(i))^2]', \qquad (8.18)$$

$$2\sum_{i=0}^{N-1} y(i) - \sum_{i=0}^{N-1}[(y'(i))^2]' \xrightarrow{Z_m} 2\sum_{k=0}^{N-1} Y(k) - \sum_{i=0}^{N-1}[(Y'(i))^2]'. \qquad (8.19)$$

Из выражения (8.19) вытекает уникальное равенство, которое не может быть получено никаким другим путем, только через двойственность:

$$2[\sum_{i=0}^{N-1} y(i) - \sum_{k=0}^{N-1} Y(k)] = \sum_{i=0}^{N-1}[(y'(i))^2]' - \sum_{i=0}^{N-1}[(Y'(i))^2]'. \qquad (8.20)$$

Из этого равенства следует, что удвоенная разность средних значений оригинала изображения равна разности средних значений от квадрата обратных значений оригинала и изображения соответственно.

В заключение следует отметить, что каждая из теорем **S**-преобразования имеет двойственность, так как все из них, так или иначе, содержат, по крайней мере, операцию сложения или умножения. Вместе с тем описанная двойственность операций, заданных на структуре (7.2), подтверждает тот факт, что доказательство теорем **S**-преобразования через представления конечных групп и свойства характеров этих групп может быть некорректным.

### 9. Теоретико-числовое преобразование над кортежами.

Кортежем называется конечная последовательность, допускающая повторение, элементов какого-либо множества. Самым простым кортежем есть упорядоченная пара, обобщением понятия которой и является понятие кортежа. Иногда упорядоченную пару называют комплексом. Кортеж трех объектов называют упорядоченной тройкой или триплексом и т. д.

И так, упорядоченной парой или просто парой называют любой кортеж, состоящий из двух элементов

$$\ddot{A} = <a, b> \qquad (9.1)$$

где $\ddot{A}$ -символ, обозначающий пару,

$a$, $b$- числа, составляющие эту пару.

Если же составить множество **P** из таких элементов, то над ними можно задать следующие операции:

сложение (вычитание) пары $\quad \ddot{A} \pm \ddot{B} = <a, b> \pm <c, d> = <a \pm c, b \pm d>;$ \qquad (9.2)

умножение пары $\quad \ddot{A} \times \ddot{B} = <a, b> \times <c, d> = <a \times c, b \times d>;$ \qquad (9.3)

деление пары $\quad \ddot{A} / \ddot{B} = <a, b>/<c, d> = <a/c, b/d>;$ \qquad (9.4)

кросс-произведение (свертка) $\quad \ddot{A} * \ddot{B} => \begin{bmatrix} a & b \\ b & a \end{bmatrix} \times \begin{bmatrix} c \\ d \end{bmatrix} = \begin{bmatrix} ac+bd \\ bc+ad \end{bmatrix}.$ \qquad (9.5)

Ясно, что операции (9.3), (9.4) не представляют существенного теоретического интереса, так как являются простым механическим повторение двух операций умножения, хотя и на практике очень часто встречаются. Другое дело кросс-произведение (9.5). Однако матрица в этом выражении такова, что ее определитель равен нулю при *a=b*. Следовательно, операция (9.5) не обратима. С другой стороны при таком умножении произведение (9.5) будет представлять пару, которая на плоскости будет принадлежать только одному квадранту. Поэтому на практике в матрице изменяют знак плюс одного из элементов матрицы на минус. В этом случае существует четыре варианта.

$$
\begin{aligned}
1.\ \ddot{A}*\ddot{B} &=> \begin{bmatrix} a & -b \\ b & a \end{bmatrix} \times \begin{bmatrix} c \\ d \end{bmatrix} = \begin{bmatrix} ac-bd \\ bc+ad \end{bmatrix}. \\
2.\ \ddot{A}*\ddot{B} &=> \begin{bmatrix} a & b \\ b & -a \end{bmatrix} \times \begin{bmatrix} c \\ d \end{bmatrix} = \begin{bmatrix} ac+bd \\ bc-ad \end{bmatrix}. \\
3.\ \ddot{A}*\ddot{B} &=> \begin{bmatrix} -a & b \\ b & a \end{bmatrix} \times \begin{bmatrix} c \\ d \end{bmatrix} = \begin{bmatrix} -ac+bd \\ bc+ad \end{bmatrix}. \\
4.\ \ddot{A}*\ddot{B} &=> \begin{bmatrix} a & b \\ -b & a \end{bmatrix} \times \begin{bmatrix} c \\ d \end{bmatrix} = \begin{bmatrix} ac+bd \\ -bc+ad \end{bmatrix}.
\end{aligned}
\quad (9.6)
$$

Первый вариант представляет альтернативу комплексным числам, ибо произведение комплексных чисел в матричной форме записывается так

$$\dot{A}*\dot{B} = \begin{bmatrix} Re\dot{A} & j\,Im\dot{A} \\ j\,Im\dot{A} & Re\dot{A} \end{bmatrix} \cdot \begin{bmatrix} Re\dot{B} \\ j\,Im\dot{B} \end{bmatrix} = \begin{bmatrix} Re\dot{A}\cdot Re\dot{B} - Im\dot{A}\cdot Im\dot{B} \\ j(Re\dot{A}\cdot Im\dot{B} + Re\dot{B}\cdot Im\dot{A}) \end{bmatrix},$$

или без мнимой единицы

$$\dot{A}*\dot{B} => \begin{bmatrix} Re\dot{A} & -Im\dot{A} \\ Im\dot{A} & Re\dot{A} \end{bmatrix} \cdot \begin{bmatrix} Re\dot{B} \\ Im\dot{B} \end{bmatrix} = \begin{bmatrix} Re\dot{A}\cdot Re\dot{B} - Im\dot{A}\cdot Im\dot{B} \\ Re\dot{A}\cdot Im\dot{B} + Re\dot{B}\cdot Im\dot{A} \end{bmatrix},$$

откуда видно, что это выражение повторяет вариант 1.

Рассмотрим теперь вариант 2. Для получения теоретико-числового преобразования необходимо в качестве переменной *s* в выражениях (3.1) и (3.2) упорядоченную пару (3.8.1). При этом, в качестве операции сложения будем использовать операцию (9.2), в качестве операции умножения выражение (3.8.6), в качестве единицы пару вида <1,0>, а в качестве нуля пару <0,0>.

Таким образом, имеем все необходимые условия, чтобы говорить о кольце упорядоченной пары. Однако теоретико-числовое преобразование существует в кольце вычетов, для чего определим модуль кольца вычетов упорядоченной пары, подставив в выражение (3.1.12) упорядоченную пару и будем иметь

$$\ddot{M} = \sum_{m=0}^{N-1} \ddot{s}^m, \qquad (9.7)$$

$$\text{где } \ddot{s} = <a,b>, \qquad (8.8)$$

$\ddot{s}^m = (<a,b>)^m$ - степенная последовательность упорядоченных пар или кортежей,

предполагающая что $\ddot{s}^0 = <1,0>$, $\ddot{s}^1 = \begin{bmatrix} a & b \\ b & -a \end{bmatrix} \times \begin{bmatrix} 1 \\ 0 \end{bmatrix}$, $\ddot{s}^2 = \begin{bmatrix} a & b \\ b & -a \end{bmatrix} \times \begin{bmatrix} a & b \\ b & -a \end{bmatrix} \times \begin{bmatrix} 1 \\ 0 \end{bmatrix}$ и т.д.

Можно воспользоваться и другой записью в следующем виде

$$\ddot{s}^m = \left( \begin{bmatrix} a & b \\ b & -a \end{bmatrix} \right)^{m+0},$$

полагающей возведение в степень заданной матрицы при условии, что

$$\left( \begin{bmatrix} a & b \\ b & -a \end{bmatrix} \right)^0 = \begin{bmatrix} 1 & 0 \\ 0 & 0 \end{bmatrix}.$$

Тогда пара преобразований запишется

$$\ddot{X}(k) \stackrel{\ddot{\mathcal{Z}}_m}{=\!=} \sum_{i=0}^{N-1} \ddot{x}(i) \ddot{s}^{-(ki) \bmod N}, \qquad (9.9)$$

$$\ddot{x}(i) \stackrel{\ddot{\mathcal{Z}}_m}{=\!=} \frac{1}{N} \sum_{k=0}^{N-1} \ddot{X}(k) \ddot{s}^{(ki) \bmod N}, \qquad (9.10)$$

$$\ddot{M} = \sum_{m=0}^{N-1} \ddot{s}^m,$$

где $\ddot{\mathcal{Z}}_m$ - кольцо вычетов упорядоченной пары по модулю,

$N$ - некоторое простое число из множества **N**, представляющее размерность преобразования,

$\stackrel{\ddot{\mathcal{Z}}_m}{=\!=}$ - означает равенство (сравнение) в кольце вычетов $\ddot{\mathcal{Z}}_m$,

$\ddot{x}(i), \ddot{X}(k)$ - последовательности упорядоченных пар, представляющие оригинал и изображение соответственно,

$i, k$ - номера (индексы) компонент последовательностей.

Выражение (9.9) представляет собой прямое преобразование, а (9.10) - обратное.

Покажем, что эта пара выражений действительно представляет собой теоретико-числовое преобразование. Заменив $i$ на $n$, во избежании путаницы индексов, и подставив его в (3.8.9), будем иметь

$$\ddot{x}(n) \stackrel{\ddot{\mathcal{Z}}_m}{=\!=} \frac{1}{N} \sum_{k=0}^{N-1} (\sum_{i=0}^{N-1} \ddot{x}(i) s^{-(ki) \bmod N}) \ddot{s}^{(kn) \bmod N}.$$

Изменив порядок суммирования, получим

$$\ddot{x}(n) \stackrel{\ddot{\mathcal{Z}}_m}{=\!=} \frac{1}{N} \sum_{k=0}^{N-1} \ddot{x}(i) (\sum_{i=0}^{N-1} \ddot{x}(i) \ddot{s}^{(kn-ki) \bmod N}).$$

Далее $$\ddot{x}(n) \stackrel{\ddot{\mathcal{Z}}_m}{=\!=} \frac{1}{N} \sum_{i=0}^{N-1} \ddot{x}(i) (\frac{1}{N} \sum_{k=0}^{N-1} \ddot{s}^{|(n-i)k| \bmod N}). \qquad (9.11)$$

Рассмотрим теперь внутреннюю сумму $$\sum_{k=0}^{N-1} \ddot{s}^{|(n-i)k| \bmod N}. \qquad (9.12)$$

Не смотря то, что переменная $\ddot{s}$ есть упорядоченная пара, при $i = n$ значение суммы

(9.12) будет равно $N$. А для $i \neq n$ внутренняя сумма (9.12) не равна $N$ и не равна нулю (что необходимо для преобразования). Однако, на основании теоремы (3.1), существует сравнение вида
$$\sum_{k=0}^{N-1} \dddot{s}^{(n-i)k | \text{mod } N} = 0 \, | \, \text{mod} \, \dddot{M}$$

в силу тех обстоятельств, что среднее значение весовой функции за период в кольце вычетов равно нулю. Тогда выражение (9.11) можно переписать в виде системы следующим образом

$$\dddot{x}(n) \stackrel{\dddot{\mathcal{Z}}_m}{=\!=\!=} \frac{1}{N} \sum_{i=0}^{N-1} \dddot{x}(i) = N, \quad i = n$$

$$\dddot{x}(n) \stackrel{\dddot{\mathcal{Z}}_m}{=\!=\!=} \frac{1}{N} \sum_{i=0}^{N-1} \dddot{x}(i) = 0, \quad i \neq n$$

Отсюда следует, что и правая часть выражения (9.11) будет состоять из единственного, не равного нулю в кольце вычетов $\dddot{\mathcal{Z}}_m$, члена $\dddot{x}(i)$ только в том случае, если $i = n$. Тогда в (9.11) последовательность $\dddot{x}(n)$ совпадает с последовательностью $\dddot{x}(i)$ только при $i = n$, следовательно, выражения (9.9) и (9.10) представляют собой теоретико-числовое преобразование.

Аналогичным образом можно показать, что существуют теоретико-числовые преобразования над упорядоченными парами и для вариантов умножения пар 3 и 4.

Литература.

1. Рабинер Л., Гоулд Б. Теория и применение цифровой обработки сигналов. - М.: Мир, 1978, - 848 с.

2. Nussbauwer H.I. Fast Fourier Transform fnd Convolution Algorithms.- Berlin, Heedelberg, New York.: Springer-Verlad, - 1981.-250 с.

3. Вариченко Л.В., Лабунец В.Г., Раков М.А. Абстрактные алгебраические системы и цифровая обработка сигналов.- Киев.: Наук. Думка, 1986. - 247 с.

4. Семотюк М.В., Палагин А.В. Обощение теоретико-числовых преобразований.//Комп`ютерні засоби, мережі та системи.- Киев: ИК НАНУ, 2002, №1. – стр. 3-13.

5. Акушский И.Н., Юдицкий Д.И. Машинная арифметика в остаточных классах. - М.: Сов. Радио, 1968. - 439 с.

6. Семотюк М.В. Обобщенное теоретико-числовое преобразование. - Киев: ИК НАН Украины, 1994. - 30 с. (Препринт/НАН Украины, Ин-т кибернетики им.В.М.Глушкова,: 94-6).